\documentclass[]{amsart}
\usepackage{amssymb}
\usepackage{latexsym}
\newtheorem{prop}{Proposition}
\newtheorem{thm}{Theorem}
\newtheorem{cor}{Corollary}
\newtheorem{lemma}{Lemma}
\newtheorem{conj}{Conjecture}
\theoremstyle{remark}

\def\P{{\mathbb P}}

\def\Z{{\mathbb Z}}

\def\R{{\mathbb R}}
\def\C{{\mathbb C}}

\def\O{{\mathcal O}}

\def\H{{\mathcal H}}
\def\l{{\lambda}}

\def\a{{\alpha}}

\def\endproof{\hfill $\Box$}

\newcommand{\dis}{\displaystyle}
\newcommand{\gequ}{\geqslant}
\newcommand{\lequ}{\leqslant}
\newcommand{\lra}{\longrightarrow}
\newcommand{\ra}{\rightarrow}

\newcommand{\Spec}{\mbox{Spec}}
  
\newcommand{\Ker}{\mbox{Ker}}
\newcommand{\ov}{\overline}
\newcommand{\noin}{\noindent}
\newcommand{\wt}{\widetilde}
\newcommand{\wh}{\widehat}

\newcommand{\med}{\medskip}

\begin{document}

\title[Standard Conjectures for the Arithmetic Grassmannian]
{Standard Conjectures for the Arithmetic Grassmannian $G(2,N)$
and Racah Polynomials}
\author{Andrew Kresch and Harry Tamvakis}
\date{March 28, 2000}
\address{Department of Mathematics, University of Pennsylvania,
209 South 33rd Street, 
Philadelphia, PA 19104, USA}
\email{kresch@math.upenn.edu}
\address{Department of Mathematics, University of Pennsylvania,
209 South 33rd Street, 
Philadelphia, PA 19104, USA}
\email{harryt@math.upenn.edu}

\maketitle 

\setcounter{section}{-1}

\section{Introduction}
\noindent

Let $X$ be an arithmetic variety, by which we mean a regular,
projective and flat scheme over $\Spec\,\Z$, of absolute 
dimension $d+1$.
Assume that $\ov{M}=(M,\|\cdot\|)$ is a hermitian line bundle
on $X$ which is arithmetically ample, in the sense of \cite{Z} 
and \cite[\S 5.2]{S}. For each $p$ the line
bundle $\ov{M}$ defines an arithmetic Lefschetz operator
\[
\begin{array}{cccc}
\wh{L}\,: & \wh{CH}^p(X)_{\R} & \lra &  \wh{CH}^{p+1}(X)_{\R} \\
    & \a  & \longmapsto   & \a \cdot \wh{c}_1(\ov{M}).
\end{array}
\]
Here $\wh{CH}^*(X)_{\R}$ is the real arithmetic Chow ring of 
\cite{GS} and $\wh{c}_1(\ov{M})$ is the arithmetic first Chern 
class of $\ov{M}$. 

In this setting,
Gillet and Soul\'e \cite{GS} proposed arithmetic analogues
of Gro\-then\-dieck's standard conjectures \cite{Gr} on algebraic 
cycles.
A more precise version of the conjectures was formulated in
\cite[\S 5.3]{S}; assuming $2p\lequ d+1$, the statement is

\begin{conj}
\label{conj1}
{\em (a)} (Hard Lefschetz) The map
\[
\wh{L}^{d+1-2p}:\wh{CH}^p(X)_{\R}\lra\wh{CH}^{d+1-p}(X)_{\R}
\]
is an isomorphism;

\noin
{\em (b)} (Hodge index) If the nonzero $x\in \wh{CH}^p(X)_{\R}$
satisfies $\wh{L}^{d+2-2p}(x)=0$, then
\[
(-1)^p\,\wh{\deg}(x\,\wh{L}^{d+1-2p}(x))>0.
\]
\end{conj}

We study these conjectures when $X=G(r,N)$ is the arithmetic
Grassmannian, parametrizing $r$-dimensional subspaces of an 
$(r+N)$-dimensional
vector space, over any field, and $\ov{M}=\ov{\O}(1)$ is the very 
ample line bundle giving the Pl\"ucker embedding, equipped with its 
natural hermitian metric. The latter is the metric induced from the 
standard metric on complex affine space, so that the first Chern form
$c_1(\ov{M})$ on $X(\C)$ is dual to the hyperplane class.

Our main result is that Conjecture \ref{conj1} holds when $r=2$.
For projective space ($r=1$) this was shown by K\"unnemann
\cite{K}. Moreover, it is proved in \cite{KM} and \cite{Ta} that 
Conjecture \ref{conj1} holds for $G(r,N)$ after a suitable scaling 
of the metric on $\ov{\O}(1)$. To obtain the precise 
result for $G(2,N)$ we use the
arithmetic Schubert calculus of \cite{T} and linear algebra 
to reduce the problem to 
combinatorial estimates. In this case the inequality in part (b)
asserts the positivity of a linear combination of harmonic numbers
with coefficients certain Racah polynomials. The latter are a system
of orthogonal polynomials in a discrete variable introduced by Wilson
\cite{Wi}\cite{AW} which generalize the classical Racah
coefficients or $6$-$j$ symbols \cite{R} of quantum physics.

The results of K\"unnemann \cite{K} show that each statement in
Conjecture \ref{conj1} (for given $X$, $p$ and $\ov{M}$) is true  
if and only if it holds when
$\wh{CH}^p(X)_{\R}$ is replaced by the Arakelov subgroup 
$CH^p(\ov{X})_{\R}$ associated to the K\"ahler form $c_1(\ov{M})$.
We therefore restrict attention to this subgroup throughout the paper.
In Section \ref{ccs} we study arithmetic Lefschetz theory for 
varieties which admit a cellular decomposition and derive a 
cohomological criterion (Corollary \ref{hlhicor}) which we use
to check Conjecture \ref{conj1}. This criterion
does not suffice to check the Hodge index inequality on more
general Grassmannians. In Section \ref{ag2n} we
apply classical and arithmetic Schubert calculus to reduce the
conjecture for $G(2,N)$ to estimates for a class of
Racah polynomials.
The required bounds for these polynomials are established in
Section \ref{bfrp}. 

We thank Christophe Soul\'e for suggesting this problem to us and for
general encouragement. Thanks are also due to Klaus K\"unnemann,
Jennifer Morse and Herb Wilf for helpful discussions.
It is a pleasure to acknowledge the support of 
National Science Foundation Postdoctoral Research Fellowships for 
both authors.

\section{Arithmetic standard conjectures on cellular spaces}
\label{ccs}

We study Conjecture \ref{conj1} for arithmetic varieties $X$
which have a cellular decomposition over $\Spec\,\Z$, in the sense of 
\cite[Ex.\ 1.9.1]{F}; the Grassmannian $G(r,N)$ is a typical example.
See \cite{KM} for more information on these 
spaces and an approach to a weaker version of the conjecture.
Recall that for each $p$ the class map
\[
\mathrm{cl}:\,CH^p(X)_{\R}\lra H^{p,p}(X_{\R})
\]
is an isomorphism of the real Chow ring 
$CH^p(X)_{\R}=CH^p(X)\otimes_{\Z}\R$ with the 
space $H^{p,p}(X_{\R})$ of real harmonic differential 
$(p,p)$-forms on $X(\C)$. We denote by
\[
\begin{array}{cccc}
L\,: & CH^p(X)_{\R} & \lra &  CH^{p+1}(X)_{\R} \\
    & \a  & \longmapsto   & \a \cdot c_1(M)
\end{array}
\]
the classical Lefschetz operator associated to an ample 
line bundle $M$ over $X$.

Let us equip the holomorphic line bundle $M(\C)$ with a smooth 
hermitian metric, invariant under complex
conjugation, to obtain a hermitian line bundle $\ov{M}$. 
As we have indicated, to check Conjecture \ref{conj1} for the
operator $\wh{L}(\alpha)=\alpha\cdot \wh{c}_1(\ov{M})$
it suffices to work with the Arakelov Chow group $CH^p(\ov{X})_{\R}$ 
defined using the K\"ahler form $c_1(\ov{M})$.
Since $X$ has a cellular decomposition, we have an exact sequence
\begin{equation}
\label{exsq}
0 \lra CH^{p-1}(X)_{\R} \stackrel{\wt{a}}\lra CH^p(\ov{X})_{\R}
\stackrel{\zeta}\lra CH^p(X)_{\R}\lra 0
\end{equation}
(see \cite[Prop.\ 6]{KM}).
Here $\wt{a}=a\circ\mathrm{cl}$ is the composite of the class
map with the
natural inclusion $a: H^{p-1,p-1}(X_{\R})\hookrightarrow CH^p(\ov{X})_{\R}$
and $\zeta$ is the projection defined in \cite[\S 1]{GS}.
We choose a splitting
\[
s_p\colon  CH^p(X)_{\R} \lra  CH^p(\ov{X})_{\R}
\]
for the sequence (\ref{exsq}) 
and thus arrive at a direct 
sum decomposition
\begin{equation}
\label{directsum}
CH^p(\ov{X})_{\R} \cong CH^p(X)_{\R}\oplus CH^{p-1}(X)_{\R}.
\end{equation}
for every $p$.

Summing (\ref{exsq}) over all $p$ produces a sequence
\begin{equation}
\label{exsq2}
0 \lra CH^{*-1}(X)_{\R} \stackrel{\wt{a}}\lra CH^*(\ov{X})_{\R}
\stackrel{\zeta}\lra CH^*(X)_{\R}\lra 0
\end{equation}
which is compatible with the actions of $L$ and $\wh{L}$.
The splitting $s:=\oplus_p s_p$ of (\ref{exsq2}) 
does not commute with $\wh{L}$ in general.
Rather, the image of $\wh{L}\circ s - s\circ L$ is 
contained in $\Ker(\zeta)$, hence
\begin{equation}
\label{defineu}
\wh{L}\circ s - s\circ L = \wt{a}\circ U
\end{equation}
for a uniquely defined degree-preserving linear operator
$U$ on $CH^*(X)_{\R}$.

We now give some conditions equivalent to the arithmetic
hard Lefschetz theorem (Theorem \ref{thmone}).
When checking these for $G(2,N)$,
we obtain something stronger, which establishes the
arithmetic Hodge index theorem as well; this is quantified in 
Theorem \ref{thmtwo}. Recall the classical Lefschetz decomposition on 
$CH^m(X)_\R\simeq H^{2m}(X(\C),\R)$:
\[
CH^m(X)_\R=\bigoplus_{p\gequ 0} L^{m-p}CH_{prim}^p(X)_\R,
\]
where the group of primitive codimension $p$ classes is
\[
CH_{prim}^p(X)_\R=\Ker(L^{d+1-2p}:CH^p(X)_{\R}\ra CH^{d+1-p}(X)_{\R}).
\]
For $m=d-p$ this decomposition induces a projection map
\[
\pi_p: CH^{d-p}(X)_{\R} \lra L^{d-2p}CH_{prim}^p(X)_{\R}.
\]

\begin{thm} \label{thmone}
Let $X$ be an arithmetic variety of dimension $d+1$
which admits a cellular decomposition.
Let $\wh L$ be the arithmetic Lefschetz operator associated to
an ample hermitian line bundle on $X$. Then the
following statements are equivalent:
\begin{itemize}
\item[(i)]
$\wh{L}^{d+1-2p}:CH^p(\ov{X})_{\R}\lra CH^{d+1-p}(\ov{X})_{\R}$
is an isomorphism for all $p$.
\item[(ii)]
There exists a linear map
$\wh{\Lambda}\colon
CH^*(\ov{X})_{\R} \lra  CH^{*-1}(\ov{X})_{\R}$
such that for every $p$ and $\alpha\in CH^p(\ov{X})_{\R}$
we have $[\wh{\Lambda},\wh{L}]\alpha=(d+1-2p)\alpha$.
\item[(iii)]
For some (equivalently, any) choice of splitting $s$ of 
{\em (\ref{exsq2})}, with $U$ as in {\em (\ref{defineu})},
\[
\pi_p\sum_{i=0}^{d-2p} L^{d-2p-i}UL^i:\,
CH_{prim}^p(X)_{\R}\lra L^{d-2p}CH_{prim}^p(X)_{\R}
\]
is an isomorphism for all $p$.
\end{itemize}
\end{thm}

\begin{proof}
We show that (iii) implies (ii).
To do this, we first prove that if
$s'$ is another splitting of (\ref{exsq2}) as above,
with associated linear operator $U'$ on $CH^*(X)_{\R}$,
then for any $p$,
\begin{equation}
\label{pfaa}
\pi_p\sum_{i=0}^{d-2p} L^{d-2p-i}U'L^i(\alpha)=
\pi_p\sum_{i=0}^{d-2p} L^{d-2p-i}UL^i(\alpha)
\end{equation}
for all $\alpha\in CH^p_{prim}(X)_{\R}$.
Indeed, from (\ref{defineu}) we have
\begin{equation}
\label{pfab}
\wh{L}^k\circ s - s\circ L^k = \wt{a}
\sum_{i=0}^{k-1} L^{k-1-i}UL^i
\end{equation}
for all $k$.
Taking $k=d-2p+1$ and using the fact that $\alpha$ is primitive,
(\ref{pfab}) gives
\begin{equation}
\label{pfac}
\wh{L}^{d-2p+1}(s'(\alpha)-s(\alpha))=
\wt{a}\sum_{i=0}^{d-2p} L^{d-2p-i}(U'-U)L^i(\alpha).
\end{equation}
But $s'(\alpha)-s(\alpha)$ is the class of a pure harmonic form,
so the left-hand side of (\ref{pfac}) is the harmonic form
associated to an element which is in the image of $L^{d-2p+1}$,
and hence is killed by $\pi_p$.

We now claim there exists a splitting $s'$ such that
for any $p$ and $\alpha\in CH^p_{prim}(X)_{\R}$ we have
\begin{equation}
\label{conda}
U'L^i\alpha = 0\, \text{ for all }\, i<d-2p \quad\text{and}\quad
U'L^{d-2p}\alpha \in L^{d-2p}CH^p_{prim}(X)_{\R}.
\end{equation}
Indeed, if we let $D$ be the linear transformation
such that $s'-s=\wt{a}\circ D$, then
$$U'=U+[L,D],$$
and it is an exercise to check that the space of
transformations $[L,D]$ is equal to the set of
operators $V$ on $CH^*(X)_{\R}$ satisfying
$\pi_p\sum_{i=0}^{d-2p} L^{d-2p-i}VL^i(\alpha)=0$
for all $\alpha\in CH^p_{prim}(X)_{\R}$ and every $p$.

Suppose (iii) holds and choose 
$s'$ satisfying (\ref{conda}).
Let us choose a primitive basis for $CH^*(X)_{\R}$.
Applying $s'$, we get half of a basis for $CH^*(\ov{X})_{\R}$.
By (iii), we may apply
$(L^{d-2p})^{-1}\pi_pU'L^{d-2p}$
to the basis elements in $CH_{prim}^p(X)_{\R}$ for each $p$
to obtain another basis for $CH^*(X)_{\R}$, which we view
(via $\wt{a}$) as the other half of our basis for
$CH^*(\ov{X})_{\R}$.

Let $v\in CH_{prim}^p(X)_{\R}$ be one of the basis elements,
and let $r=d-2p$.
By our conditions on $s'$,
a subset of our basis for $CH^*(\ov{X})_{\R}$ consists of
$\wh{v}:=s'(v)$, the iterates $\wh{L}^i(\wh{v})=s'(L^iv)$
of $\widehat L$ applied to $\wh{v}$,
the primitive element $w$ satisfying
$L^r(w)=\pi_p(U'L^r(v))$,
and the iterates of $\widehat L$ applied to $w$:
\begin{equation}
\label{hjry}
\wh{v},\, \wh{L}\wh{v},\, \ldots,\, 
\wh{L}^r\wh{v},\, w,\, Lw,\, \ldots,\, L^rw.
\end{equation}
The action of $\wh{L}$ is to send each element in (\ref{hjry})
to the element on its right, except that 
$\wh{L}^r\wh{v}$ is sent to $L^rw$,
and $L^rw$ to $0$.
We construct $\wh{\Lambda}$ explicitly: define
\begin{align*}
\wh{\Lambda}(\wh{L}^i\wh{v})&=i(r+2-i)\wh{L}^{i-1}\wh{v},\\
\wh{\Lambda}(L^iw)&=(r+1)\wh{L}^i\wh{v} + i(r-i)L^{i-1}w.
\end{align*}
Then $\wh{\Lambda}$ (defined this way for every basis
element $v$) satisfies the condition of (ii).

Statement (i) follows from (ii) by standard representation theory of 
$sl(2)$, as in \cite[\S V.3]{W}.
To show that (i) implies (iii), note that by (\ref{pfaa}) the
condition in (iii) is independent of choice of splitting.
If $\alpha\in \Ker \bigl(\pi_p{\textstyle\sum_{i=0}^{d-2p}
L^{d-2p-i}UL^i\bigr)}$ is a nonzero primitive element 
and if we take $s'$ to be a splitting which satisfies (\ref{conda}),
then
$$\wh{L}^{d+1-2p}(s'(\alpha))=\bigl(\wh{L}^{d+1-2p}\circ s' - 
s'\circ L^{d+1-2p}
\bigr) \alpha = \sum_{i=0}^{d-2p} L^{d-2p-i}U'L^i(\alpha)
=0,$$
and (i) fails.
\end{proof}

\begin{thm}
\label{thmtwo}
Suppose the arithmetic variety $X$ and $p$
are such that $CH^{p-1}_{prim}(X)_{\R}=0$.
If, for each
nonzero $\alpha\in CH^p_{prim}(X)_{\R}$, we have
\begin{equation}
\label{needforhi}
(-1)^p \sum_{i=0}^{d-2p} \int_X L^{d-2p-i}\alpha\wedge UL^i\alpha > 0
\end{equation}
then the statements in the arithmetic hard Lefschetz and Hodge index
conjectures are true for that $X$, $p$ and $\ov{M}$.
\end{thm}

\begin{proof}
Let
$(\alpha,\beta)\in CH^p(\ov{X})_{\R}$ be a nonzero element of 
the kernel of $\wh{L}^{d+2-2p}$;
the notation $(\alpha,\beta)$ refers to the direct sum decomposition
(\ref{directsum}), with respect to some splitting. We claim that
$\alpha$ must be in $CH^p_{prim}(X)_{\R}$.
Indeed, $\wh{L}^{d+2-2p}(\alpha,\beta)=(L^{d+2-2p}\alpha,\gamma)$
for some $\gamma$ and 
$L^{d+2-2p}\alpha=0$ implies $L^{d+1-2p}\alpha=0$ since
$CH^{p-1}_{prim}(X)_{\R}$ vanishes.
Also, by the classical hard Lefschetz theorem, $\alpha\ne 0$.
Now, if 
\[
\langle\,\ ,\ \rangle\, :\, CH^*(\ov{X})_{\R}\otimes 
CH^*(\ov{X})_{\R}\lra \R
\]
denotes the arithmetic intersection pairing, then we have
\begin{align*}
\langle\,(\alpha,\beta)\, , \,\wh{L}^{d+1-2p}(\alpha,\beta) \,\rangle 
&= \langle\,(\alpha,\beta)\, , (0,
\sum_i L^{d-2p-i}UL^i\alpha + L^{d+1-2p}\beta)\, \rangle \\
&= \frac{1}{2}\sum_i \int_X L^{d-2p-i}\alpha\wedge UL^i\alpha.
\end{align*}
Hence, assuming $CH^{p-1}_{prim}(X)_{\R}=0$,
we have $\wh{L}^{d+1-2p}(\alpha,\beta)\ne 0$ for every nonzero
$(\alpha,\beta)\in CH^p(\ov{X})_{\R}$.
Moreover, if $(\alpha,\beta)$ is primitive,
then the pairing of $(\alpha,\beta)$ with
$\wh{L}^{d+1-2p}(\alpha,\beta)$ has the required sign.
\end{proof}

\begin{cor}
\label{hlhicor} 
Suppose $X$ is such that, for every $p$, 
\begin{equation}
\label{primcond}
CH^p_{prim}(X)_{\R}\ne 0 \qquad\text{implies}\qquad
CH^{p-1}_{prim}(X)_{\R}=0.
\end{equation}
If condition {\em (\ref{needforhi})} holds 
for every $p$ and each 
nonzero $\alpha\in CH^p_{prim}(X)_{\R}$, then 
both the arithmetic hard Lefschetz and Hodge index
conjectures are true for $X$, $\ov{M}$.
\end{cor}

\medskip
\noin
{\bf Example.} We illustrate both theorems for projective
space $\P^n$ over $\Spec\,\Z$
(compare \cite[\S 4]{K}). In this case we 
choose the splitting
\[
CH^*(\ov{\P^n})_{\R}=\bigoplus_{i=0}^n\R\cdot \wh{\omega}^i \,\oplus\,
\bigoplus_{i=0}^n\R\cdot \omega^i
\]
where $\wh{\omega}^i=\wh{c}_1(\ov{\O}(1))^i$ and 
$\omega^i=\wt{a}(c_1(\O(1))^i)$.
Then the sequence (\ref{hjry}) is given by
\[
\wh{1},\, \wh{\omega},\, \ldots,\, \wh{\omega}^n,\,
\tau_n,\, \tau_n\omega,\, \ldots,\, \tau_n\omega^n.
\]
Here $\tau_n=\sum_{k=1}^n\H_k$, where each 
$\H_k=1+\frac{1}{2}+\cdots
+\frac{1}{k}$ is a harmonic number. In the proof of Theorem
\ref{thmone} we constructed an explicit adjoint map
$\wh{\Lambda}$ for the arithmetic Lefschetz
operator $\wh{L}(x)=\wh{\omega}\cdot x$; in
our example it is given by
\begin{align*}
\wh{\Lambda}(\wh{\omega}^i)&=i(n+2-i)\,\wh{\omega}^{i-1},\\
\wh{\Lambda}(\omega^i)&=\frac{n+1}{\tau_n}\,\wh{\omega}^i + 
i(n-i)\,\omega^{i-1}.
\end{align*}
Observe that the nonzero primitive elements of $CH^*(\P^n)_{\R}$
are multiples of $1\in CH^0(\P^n)_{\R}$; hence $\P^n$ satisfies
(\ref{primcond}). The operator $U$ is given by
$U(\omega^i)=\delta_{i,n}\tau_n\omega^n$, and condition
(\ref{needforhi}) for $p=0$, $\a=1$ becomes
\[
\sum_{i=0}^n\int_{\P^n}\omega^{n-i}\wedge U(\omega^i)
=\tau_n\int_{\P^n}\omega^n=\tau_n>0.
\]
The arithmetic Hodge index conjecture for $\P^n$ follows by 
applying Corollary \ref{hlhicor}.

\section{The arithmetic Grassmannian $G(2,N)$}
\label{ag2n}
\noindent

In this section we study Conjecture \ref{conj1}
for the Grassmannian of lines in projective space.
For computational 
purposes we will work with the isomorphic Grassmannian 
$G=G(N,2)$ parametrizing $N$-planes in $(N+2)$-space throughout.
Note that $d=\dim_{\C}G(\C)=2N$. 
There is a universal exact sequence of vector bundles
\[
0 \lra S \lra E\lra Q \lra 0
\]
over $G$; the complex points of $E$ and $Q$ 
are metrized by giving
the trivial bundle $E(\C)$ the trivial hermitian
metric and the 
quotient bundle $Q(\C)$ the induced metric.
The hermitian vector bundles that result are denoted $\ov{E}$

The real vector space $CH^*(G)_{\R}
\cong H^{2*}(G(\C),\R)$ decomposes as
\[
CH^*(G)_{\R}=\bigoplus_{a,b}\R\cdot s_{a,b}(Q),
\]
summed over all partitions $\l=(a,b)$ with $a\lequ N$,  i.e., whose
Young diagrams are contained in the $2\times N$ rectangle $(N,N)$.
Moreover $s_{\l}(Q)=s_{a,b}(Q)$
is the characteristic class coming from the Schur polynomial
$s_{a,b}$ in the Chern roots of $Q$; this is dual to the class
of a codimension $|\l|=a+b$ Schubert variety in $G$.

The line bundle $M=\det(Q)$ giving the Pl\"ucker embedding has
$c_1(M)=s_1(Q)$; let $L:CH^p(G)_{\R}\ra CH^{p+1}(G)_{\R}$ be the 
associated classical Lefschetz operator. Further for all $p$ let
 $*:CH^p(G)_{\R}\ra CH^{2N-p}(G)_{\R}$ 
denote the Hodge star operator induced by the K\"ahler form 
$s_1(\ov{Q})$. We then have

\begin{prop}
\label{clprim}
The space $CH^p_{prim}(G)_{\R}$ is nonzero if and only if $p=2k\lequ
N$. In the latter case it is one dimensional and spanned by the class
\[
\alpha_k=\sum_{j=0}^k(-1)^j\binom{N+1-j}{N-2k}\binom{N-2k+j}{N-2k}
s_{2k-j,j}(Q).
\]
\end{prop}

\noin
{\em Proof.} By computing the Betti numbers for $G$ one sees that
\[
\dim CH^p(G)_{\R}-\dim CH^{p-1}(G)_{\R}>0 
\]
if and only if $p=2k\lequ N$, and in this case the above difference
equals 1. For such $p$ we have
\begin{equation}
\label{kerL}
\Ker(L:CH^{2N-2k}(G)_{\R}\ra CH^{2N-2k+1}(G)_{\R})=\mathrm{Span}\{
\sum_{j=0}^k(-1)^j s_{N-j,N-2k+j}(Q)\}.
\end{equation}
One checks (\ref{kerL}) easily using the Pieri rule:
\[
L(s_{a,b}(Q))= s_{a+1,b}(Q)+s_{a,b+1}(Q)
\]
where it is understood that $s_{c,c'}(Q)=0$ if $c<c'$ or $c>N$.

{}From \cite{KT} we know the action of the Hodge star operator on
$CH^*(G)_{\R}$ is given by
\begin{equation}
\label{star}
*\,s_{a,b}(Q) = \frac{(a+1)!\,b!}{(N-a)!\,(N-b+1)!} \, s_{N-b,N-a}(Q).
\end{equation}
Since
\[
CH^{2k}_{prim}(G)_{\R}= * \,
\Ker\left(L:CH^{2N-2k}(G)_{\R}\ra CH^{2N-2k+1}(G)_{\R}\right),
\]
the proof is completed by applying (\ref{star}) to (\ref{kerL}) and
noting that the result is proportional to $\a_k$. \endproof

\medskip

We now pass to the arithmetic setting, where we use the arithmetic
Schubert calculus of \cite[\S 3,4]{T}. The real Arakelov Chow group
$CH^p(\ov{G})_{\R}$ decomposes as
\begin{equation}
\label{areq}
CH^p(\ov{G})_{\R} = 
\bigoplus_{a+b=p}\R\cdot \wh{s}_{a,b}(\ov{Q}) \,
\oplus \bigoplus_{a'+b'=p-1}\R\cdot
s_{a',b'}(\ov{Q}).
\end{equation}
Here the indexing sets satisfy $N\gequ a \gequ b\gequ 0$,
$\wh{s}_{a,b}(\ov{Q})$
is an arithmetic characteristic class and we identify the 
harmonic differential form $s_{a',b'}(\ov{Q})$ with its image
in $CH^p(\ov{G})_{\R}$.
The decomposition (\ref{areq}) is induced by the splitting map
$s_{a,b}(Q)\longmapsto \wh{s}_{a,b}(\ov{Q})$ which agrees with
the one used in \cite{T}.

The hermitian line bundle $\ov{M}$ 
has $\wh{c}_1(\ov{M})=\wh{s}_1(\ov{Q})$ and
is arithmetically ample; this follows from \cite[Prop.\ 3.2.4]{BGS}.
We now apply the arithmetic Pieri rule of \cite[\S 4]{T} to 
compute the action of the arithmetic Lefschetz operator 
$\wh{L}(x)=\wh{s}_1(\ov{Q})\cdot x$ on the above basis elements.
The induced operator $U:CH^*(G)_{\R}\ra CH^*(G)_{\R}$ 
of (\ref{defineu}) satisfies $U(s_{a,b})=0$ for $a<N$ and 
\begin{equation}
\label{Ueq}
U(s_{N,b})=\Bigl(\sum_{i=0}^{N+1}\H_i\Bigr)s_{N,b} -
\sum_{i=0}^{\lfloor (N-b)/2\rfloor} (\H_{N-b+1-i} - \H_i ) s_{N-i, b+i}.
\end{equation}
Here and in the rest of this section
$s_{a,b}$ will denote the Schubert class
$s_{a,b}(Q)\in CH^{a+b}(G)_{\R}$ and $\H_i$ is a harmonic number;
by convention $\H_0=0$. Recall that the classical intersection
pairing on $CH(G)_{\R}$ satisfies
\[
\left<s_{a,b},s_{a',b'}\right>=\int_Gs_{a,b}\wedge s_{a',b'}
=\delta_{(a,b),(N-b',N-a')}.
\]

The sequence of Betti numbers for $G$ shows that $G$ satisfies 
condition (\ref{primcond}) of Corollary \ref{hlhicor}. 
We proceed to check the inequality (\ref{needforhi}) for all even 
$p=2k$; this will establish Conjecture \ref{conj1} for $G(N,2)$. In our
case (\ref{needforhi}) may be written as
\[
\Sigma(N,k):=\int_G\sum_{b=0}^{N-2k}
L^{N-2k-b}\a_k\wedge UL^{N-2k+b}\a_k>0.
\]

To compute
iterates of the classical Lefschetz operator $L$ on the Schubert
basis, note that
\begin{equation}
\label{Liter}
L^rs_{\mu}=\sum_{{\l\supset \mu} \atop {|\l|=|\mu|+r}}
f^{\l/\mu}s_{\l}.
\end{equation}
When $\l=(\l_1,\l_2)$ and $\mu=(\mu_1,\mu_2)$ are partitions with
at most two parts
the skew $f$-number in (\ref{Liter}) satisfies
\begin{equation}
\label{bindiff}
f^{\l/\mu}={|\l|-|\mu| \choose \l_1-\mu_1} - 
{|\l|-|\mu| \choose \l_1-\mu_2+1}.
\end{equation}
This follows from the determinantal formula for $f^{\l/\mu}$, given
for example in \cite[Corollary 7.16.3]{St}.

We now apply Proposition \ref{clprim} and (\ref{Liter}), (\ref{bindiff})
to calculate
\begin{gather*}
L^{c-2k}\a_k = \sum_{j=0}^k(-1)^j
\binom{N+1-j}{N-2k}\binom{N-2k+j}{N-2k}
L^{c-2k}s_{2k-j,j} \\
= \sum_{j=0}^k\sum_{i\gequ j}(-1)^j\binom{N+1-j}{2k+1-j}\binom{N-2k+j}{j}
\left[\binom{c-2k}{i-j}-\binom{c-2k}{i-2k-1+j}\right]s_{c-i,i}.
\end{gather*}
Therefore,
\begin{equation}
\label{zeroeq}
L^{c-2k}\a_k = \sum_{i,j}(-1)^j\binom{N+1-j}{N-2k}\binom{N-2k+j}{N-2k}
\binom{c-2k}{i-j}s_{c-i,i}.
\end{equation}
We use (\ref{zeroeq}) with $c=N+b$ to identify the coefficient of
$s_{N,b}$ in the expansion of $L^{N+b-2k}\a_k$ as
\begin{align*}
\int_GL^{N+b-2k}\a_k\wedge s_{N-b} &= 
\sum_j(-1)^j\binom{N+1-j}{N-2k}\binom{N-2k+j}{j}
\binom{N+b-2k}{N-2k+j} \\
&= \binom{N-2k+b}{N-2k}\sum_j(-1)^j\binom{b}{j}\binom{N+1-j}{N-2k} \\
&= \binom{N+1-b}{2k+1}\binom{N-2k+b}{N-2k} =:C_b.
\end{align*}
It now follows from (\ref{Ueq}) that
\begin{equation}
\label{secondeq}
UL^{N-2k+b}\alpha_k =
C_b
\Bigl(\sum_{i=0}^{N+1}\H_i\Bigr) s_{N,b} 
 - \sum_{i=0}^{\lfloor (N-b)/2\rfloor}
 C_b (\H_{N-b+1-i} - \H_i ) s_{N-i, b+i}.
\end{equation}
Note also that we have the identity
\begin{equation}
\label{msum}
\sum_{b=0}^{N-2k}C_b=\binom{2N-2k+2}{N+2}.
\end{equation}
Now we substitute $c=N-b$ in (\ref{zeroeq}), pair with 
(\ref{secondeq}) and use (\ref{msum}) to sum over $b$ and
obtain
\begin{equation}
\label{endequ}
\Sigma(N,k)=
A_{N,k}\sum_{i=1}^{N+1}\H_i+\sum_{i=0}^{\lfloor (N-b)/2\rfloor}C^i_{N,k},
\end{equation}
where 
\[
A_{N,k}=\binom{N+1}{N-2k}\binom{2N-2k+2}{N+2}
\]
and 
\[ 
C^i_{N,k}=\sum_{j,b}(-1)^j
\binom{N+1-j}{N-2k}
\binom{N-2k+j}{N-2k}
\binom{N-2k-b}{i-j}
C_b(\H_i-\H_{N-b+1-i}).
\]
By dividing the expression for $C^i_{N,k}$ into two sums and 
substituting in (\ref{endequ}) one gets
\begin{equation}
\label{eeq}
\Sigma(N,k)=
A_{N,k}\sum_{i=1}^{N+1}\H_i+\sum_{i=1}^{N+1}B^i_{N,k}\H_i,
\end{equation}
where
\[
B^i_{N,k}=\sum_{j,b}(-1)^j
\binom{N+1-j}{N-2k}
\binom{N-2k+j}{N-2k}
\binom{N-2k-b}{i-j}
C_b.
\]
(Notice that when $N-b=2r-1$ is odd, there 
is a missing summand (for $i=r$)
\[
\sum_j(-1)^j
\binom{N+1-j}{N-2k}
\binom{N-2k+j}{N-2k}
\binom{2r-2k-1}{r-j}
C_b\H_r
\]
which vanishes, as can be seen by the change of variable $j\mapsto
2k+1-j$.)

At this point it is convenient to introduce the variable change
$$n=N-2k \qquad\text{and}\qquad T=N+2$$
and write equation (\ref{eeq}) in the new coordinates as
\begin{equation}
\label{lasteeq}
\Sigma(n,T)=
A_{n,T}\sum_{i=1}^{T-1}\H_i+\sum_{i=1}^{T-1}B^i_{n,T}\H_i.
\end{equation}
Observe that
\begin{align*}
B^i_{n,T} &= 
\sum_j(-1)^j
\binom{n+j}{n}
\binom{T-1-j}{n}
\sum_b
\binom{n-b}{i-j}
\binom{T-1-b}{n-b}
\binom{n+b}{n} \\
&= \sum_j(-1)^j
\binom{n+j}{n}
\binom{T-1-j}{n}
\binom{T-1-n+i-j}{i-j}
\binom{T+n}{n-i+j}.
\end{align*}
We now substitute
$r=i-j$ and write the resulting sum in 
hypergeometric notation \cite{Ro}
\cite[Chap. 3]{VK}:
\begin{align*}
(-1)^iB^i_{n,T} &= 
\sum_r(-1)^r
\binom{n+i-r}{n}
\binom{T-1+r-i}{n}
\binom{T-1-n+r}{r}
\binom{T+n}{n-r} \\
&= \binom{n+i}{n} \binom{T+n}{n} \binom{T-1-i}{n}
{}_4F_3\biggl(\matrix -n,\, -i,\, T-n,\, T-i \\
-n-i,\, T+1,\, T-n-i  \endmatrix\biggm|1\biggr).
\end{align*}
The Whipple transformation \cite[\S 10]{Wh} applied to the above ${}_4F_3$
gives
\begin{equation}
\label{connect34}
(-1)^i\frac{B^i_{n,T}}{A_{n,T}} = 
{}_4F_3\biggl(\matrix -n,\, n+1,\, -i,\, i+1 \\
1,\, 1+T,\, 1-T \endmatrix\biggm|1\biggr).
\end{equation}
The hypergeometric term
(\ref{connect34}) belongs to 
a class of orthogonal polynomials called Racah polynomials,
which are studied in the next section.

\section{Bounds for Racah polynomials}
\label{bfrp}

 The Racah coefficients \cite{R} or $6$-$j$ symbols have long been 
used by physicists as the transformation coefficients between two
different coupling schemes of three 
angular momenta; see \cite{BL} for an exposition. In mathematical
language they are the entries of a change of basis matrix for the tensor
product of three irreducible representations of $SU(2)$; the two
bases involved come from the associativity relation for this product 
(see \cite[\S 8.4]{VK}). It was recognized later by Wilson \cite{Wi}
that these coefficients are special cases of a class of
orthogonal polynomials $R_n(x;\alpha,\beta,\gamma,\delta)$, called
Racah polynomials \cite[\S8.5]{VK}:
$$R_n(s(s+\gamma+\delta+1);\alpha,\beta,\gamma,\delta)=
{}_4F_3\biggl(\matrix -n,\, n+\alpha+\beta+1,\, -s,\, s+\gamma+\delta+1\\
\alpha+1,\, \beta+\delta+1,\, \gamma+1\endmatrix\biggm|1\biggr).$$

The Racah polynomials in (\ref{connect34}) have
$\alpha=\beta=\gamma+\delta=0$, with $\gamma=T$,
a positive integer. We let
\[
R_n(s,T)=R_n(s(s+1);0,0,T,-T).
\]
Observe that $R_n(s,T)$ is symmetric in $n$ and $s$. 
The orthogonality condition (loc.\ cit.\ or \cite{AW}) reads:
\begin{equation}
\label{ortho}
\sum_{s=0}^{T-1} (2s+1)R_n(s,T)R_m(s,T)=
\frac{T^2}{(2n+1)}\delta_{nm}.
\end{equation}

The arithmetic Hodge index inequality $\Sigma(n,T)>0$ can be rephrased
using (\ref{lasteeq}) and (\ref{connect34}) as
\begin{equation}
\label{needed}
\sum_{s=1}^{T-1} (-1)^{s+1} R_n(s,T) \H_s <
\sum_{s=1}^{T-1} \H_s.
\end{equation}
We give a proof of (\ref{needed}) which does not depend on the
precise values of the harmonic numbers. Let us say that a 
sequence $\{\H_k\}_{k\gequ 1}$ of positive real numbers 
(with $\H_0=0$) is {\em concave increasing} if $\H_k=\sum_{i=1}^kh_i$
for some monotone decreasing sequence $\{h_i\}$ of positive reals.

\begin{thm}
\label{lastthm}
Let $\{\H_k\}$ be any concave increasing sequence of real numbers and
$n$, $T$ integers with $0\lequ n\lequ T-1$ and  
$T\gequ 3$. Then inequality {\em (\ref{needed})} holds.
\end{thm}

We believe that, in fact, (\ref{needed}) holds for an arbitrary
 sequence of positive
real numbers $\H_k$, that is, the arithmetic standard conjectures for
$G(2,N)$ do not depend on the relative sizes of the harmonic numbers
involved:

\begin{conj}
\label{conj2}
For any integers $n,s$ with $0\lequ n,s\lequ T-1$
we have $| R_n(s,T)|\lequ 1$.
\end{conj}

\noin
In Proposition \ref{conjpr} we 
check this conjecture for some values of $n$
near the endpoints $0$ and $T-1$. 
Computer calculations support the validity of
Conjecture \ref{conj2} for general $n$.

\begin{proof}[Proof of Theorem \ref{lastthm}]
We shall see that (\ref{ortho}) implies (\ref{needed}) except when
$T$ is exponentially large compared to $n$.
For large $T$, the Racah polynomials are close approximations
of classical orthogonal polynomials, in this case
the Legendre polynomials, and we know how to bound these.

By Cauchy's inequality, (\ref{ortho}) gives
$$\biggl(\sum_{s=0}^{T-1}
|R_n(s,T)|\,\H_s\biggr)^2\lequ
\frac{T^2}{2n+1}
\sum_{s=0}^{T-1} \frac{\H_s^2}{2s+1}.
$$
So, (\ref{needed}) holds whenever
\begin{equation}
\label{newineq}
\sum_{s=0}^{T-1} \frac{\H_s^2}{2s+1} <
(2n+1)\biggl(\frac{1}{T}\sum_{s=0}^{T-1}\H_s\biggr)^2.
\end{equation}
Since $\{\H_k\}$ is concave increasing, the average value of
$\H_0$, $\ldots$, $\H_{T-1}$ is at least $\H_{T-1}/2$.
As $\sum_{s=1}^{T-1} 2/(2s+1)\lequ\log T$,
the inequality (\ref{newineq}) holds whenever
\begin{equation}
\label{goodrange}
\log T<n+\frac{1}{2}.
\end{equation}

To analyze the case where $T$ is exponentially
large compared to $n$,
it is convenient to introduce the change of variable
\begin{equation}
\label{change}
x=s(s+1)=-1/4+T^2(1+t)/2
\end{equation}
and the rescaling
\[
p_n(t)=(-1)^n\prod_{i=1}^n \frac{T^2-i^2}{T^2} R_n(x;0,0,T,-T).
\]
Let $P_n(t)=P^{(0,0)}_n(t)$ denote the $n^{\rm th}$ Legendre
polynomial. It is known \cite[\S3.8]{NSU} that
$$p_n(t)=P_n(t) + O(1/T^2)$$
where the constant in the error term depends on both $n$ and $t$.
For our purposes, we demonstrate
\begin{lemma}
\label{1stlemma}
{\em a)} Let $n$ and $T$ be positive integers
such that $1+2n+2n^2<T^2/10$.
Then
\begin{equation}
\label{pbound}
|p_n(t)-P_n(t)|\lequ (3/2)\cdot 4^n/T^2
\end{equation}
for all $t$ with $-1\lequ t \lequ 1$.

\med
\noin
{\em b)} We have $|p_n(t)-P_n(t)|\lequ 1/10$ 
whenever $T\gequ 90$ and $n<\log T$.
\end{lemma}

\begin{proof}
We have the following recurrences (loc.\ cit.; for $P_n(t)$ this is
classical)
\begin{align}
tP_n(t) &= \frac{n+1}{2n+1} P_{n+1}(t) + \frac{n}{2n+1}P_{n-1}(t)
\label{recurP} \\
tp_n(t) &=
\frac{n+1}{2n+1} p_{n+1}(t) -
\frac{2n^2+2n+1}{2T^2}p_n(t) +
\Bigl(1-\frac{n^2}{T^2}\Bigr)^2\frac{n}{2n+1}p_{n-1}(t)
\label{recurp}
\end{align}
with initial data
\begin{equation}
\label{initial}
\begin{array}{lll}
P_0(t)=1  &;&   P_1(t)=t  \medskip \\
p_0(t)=1  &;&   p_1(t)=t+1/(2T^2).
\end{array}
\end{equation}
Subtracting (\ref{recurP}) from (\ref{recurp}) leads to a recurrence
in $p_n(t)-P_n(t)$.
Then (\ref{pbound}) follows by induction on $n$, using
the known bound $|P_n(t)|\lequ 1$ for all $n$ and all $t$ with 
$-1\lequ t\lequ 1$. The statement (b) is a corollary of (a).
\end{proof}

\begin{prop}
\label{conjpr} We have $|R_n(s,T)|\lequ 1$ when $n\lequ 3$ or 
$n=T-1$. 
\end{prop}

\begin{proof} For $n=0$ we have $R_0(s,T)=1$. 
When $n=1$, we see from (\ref{initial})
that $p_1(t)$ is a increasing linear function in $t$,
attaining minimum when $s=0$, giving $R_1(0,T)=1$, and
maximum when $s=T-1$, giving $R_1(T-1,T)=(1-T)/(1+T)<1$.
For $n=T-1$ the Pfaff-Saalsch\"utz identity \cite{Ro} 
\cite[8.3.3]{VK} 
gives
\begin{align*}
R_{T-1}(s,T) &= \sum_{j}(-1)^j\frac{T}{T+j}\binom{s}{j}\binom{s+j}{j}
\\
&= {}_3F_2\biggl(\begin{array}{ccc} -s,\, s+1,\, T \\
1,\, T+1 \end{array}\biggm|1\biggr)  \\
&= \frac{(1-T)(2-T)\cdots (s-T)}{(1+T)(2+T)\cdots (s+T)}
\end{align*}
so the inequality is clear.

When $n=2$ or $n=3$, $p_n(t)$ is a quadratic or cubic polynomial, 
and it is a calculus exercise to check that
$|R_n(s,T)|\lequ 1$ for every integer $s$ with
$0\lequ s\lequ T-1$. In fact, the integrality condition on $s$ is
required only when $n=3$, $T=4$.
%
\end{proof}

\begin{lemma} \label{boundone} {\em a)} We have $|P_n(t)|\lequ 3/4$ for 
$t\in \R$, $|t|\lequ 0.9$ and $n\gequ 2$.

\med
\noin
{\em b)}  For $T\gequ 10$, we have $|t|\lequ 0.9$ in {\em (\ref{change})} 
whenever $\sqrt{5}/10\lequ s/T\lequ 4/5$.

\noin
{\em c)} Assume $T\gequ 90$ and $n<\log T$.
Then $\dis \frac{1}{T^{2n}}\prod_{i=1}^n(T^2-i^2)>40/41$.
\end{lemma}

\begin{proof} The indicated 
bound on Legendre polynomials  is evident for $n=2$,
and for larger $n$ it follows from the inequality
\begin{equation}
\label{exer}
\sqrt{\sin\theta}\,|P_n(\cos\theta)|\,< \,
\root\of{\frac{2}{\pi n}}, \ \ \ \ 0 \lequ \theta \lequ \pi.
\end{equation}
One obtains (\ref{exer}) by using the 
transformed differential equation for
$\sqrt{\sin\theta}\,P_n(\cos\theta)$ \cite[IV.4.24.2]{Sz};
this is indicated in \cite[Chap.\ 5, Exer.\ 15--16]{Ho}. 
The proofs of (b) and (c) are routine; the inequality
$-\log(T^{-2n}\prod_{i=1}^n(T^2-i^2))\lequ (2/T^2)
\sum_{i=1}^ni^2$ may be used for the latter.
\end{proof}

\med

To complete the proof of Theorem \ref{lastthm}, assume that 
(\ref{goodrange}) fails, so that $n< \log T -1/2$. If $T\lequ 90$
then $n\lequ 3$ and (\ref{needed}) follows from 
Proposition \ref{conjpr} (note that the inequality in the proposition
is strict unless $n=0$ or $s=0$). When $T\gequ 90$
we combine Lemma \ref{1stlemma}(b) with Lemma \ref{boundone}
to deduce the inequality (\ref{needed}).
Indeed, $(40/41)|R_n(s,T)|$ is bounded by $1+1/10$ for every $s$,
and by $3/4+1/10$ over the middle half of the summation range.
By pairing terms $\H_s$ with $\H_{T-1-s}$ and using the fact
that $\H_s+\H_{T-1-s}$ is monotone increasing for 
$0\lequ s\lequ (T-1)/2$, we obtain (\ref{needed}).
\end{proof}

\end{document}